%% file: 2004-8.tex
\documentclass{gtart_h}

\input gtoutput

\lognumber{389}
\volumenumber{8}\papernumber{8}\volumeyear{2004}
\pagenumbers{311}{334}
\received{3 December 2003}
\revised{12 February 2004}
\published{14 February 2004}
\accepted{14 February 2004}
\proposed{Robion Kirby}
\seconded{John Morgan, Ronald Stern}

\usepackage{amsmath,amsfonts,amscd}

\newcommand\Rk{\mathrm{rk}}

\newcommand{\HF}{HF}

\newtheorem{theorem}{Theorem}[section]

\newtheorem{cor}[theorem]{Corollary}

\newcommand{\Q}{\mathbb{Q}}
\newcommand{\R}{\mathbb{R}}

\newcommand{\Z}{\mathbb{Z}}

\newcommand{\Zmod}[1]{\Z/{#1}\Z}

\newcommand{\cm}{\cdot}

\newcommand{\Nbd}[1]{{\mathrm{nd}}(#1)}
\newcommand{\nbd}[1]{\Nbd{#1}}

\newcommand{\ModSWfour}{\mathcal{M}}
\newcommand{\ModFlow}{\ModSWfour}

\newcommand{\SpinC}{{\mathrm{Spin}}^c}

\newcommand\Hom{\mathrm{Hom}}

\newcommand\abuts\Rightarrow
\newcommand\Sym{\mathrm{Sym}}

\newcommand\HFpRed{\HFp_{\red}}
\newcommand\HFmRed{\HFm_{\red}}
\newcommand\mix{\mathrm{mix}}

\newcommand\x{\mathbf x}
\newcommand\w{\mathbf w}

\newcommand\y{\mathbf y}

\newcommand\ModSphere{\ModFlow\left({\mathbb S}\longrightarrow 
\Sym^{g-1}(\Sigma_{1})\times \Sym^2(\Sigma_{2})\right)}
\newcommand\ModSpheres\ModSphere
\newcommand\CF{CF}

\newcommand\CFa{\widehat{CF}}

\newcommand{\red}{\mathrm{red}}

\newcommand\HFp{\HFb}

\newcommand\HFm{\HF^-}

\newcommand\HFinf{HF^\infty}

\newcommand\HFa{\widehat{HF}}
\newcommand\HFb{HF^+}

\newcommand\Mas{\mu}
\newcommand\UnparModSp{\widehat \ModSp}
\newcommand\UnparModFlow\UnparModSp
\newcommand\Mod\ModSp

\newcommand{\spinc}{\mathfrak s}

\newcommand{\spinct}{\mathfrak t}

\newcommand\ModMaps{\mathcal M}
\newcommand\ModSp\ModMaps

\newcommand\Ta{{\mathbb T}_{\alpha}}
\newcommand\Tb{{\mathbb T}_{\beta}}
\newcommand\Tc{{\mathbb T}_{\gamma}}

\newcommand\alphas{\mbox{\boldmath$\alpha$}}

\newcommand\betas{\mbox{\boldmath$\beta$}}
\newcommand\gammas{\mbox{\boldmath$\gamma$}}

\newcommand\uCF{\underline\CF}
\newcommand\uCFinf{\uCF^\infty}
\newcommand\uHF{\underline{\HF}}
\newcommand\uHFp{\underline{\HF}^+}

\newcommand\uHFpRed{\uHFpred}
\newcommand\uHFpred{\underline{\HF}^+_{\red}}
\newcommand\uCFp{\uCF^+}
\newcommand\uCFm{\uCF^-}
\newcommand\uHFm{\underline{\HF}^-}
\newcommand\uHFa{\underline{\HFa}}
\newcommand\uCFa{\underline{\CFa}}

\newcommand\uHFinf{\uHF^\infty}

\newcommand\uFp[1]{{\underline F}^+_{#1}}

\newcommand\Fm[1]{F^{-}_{#1}}

\newcommand\Fp[1]{F^{+}_{#1}}

\newcommand\Fa[1]{{\widehat F}_{#1}}
\newcommand\Finf[1]{F^{\infty}_{#1}}

\newcommand\Fmix[1]{F^{\mix}_{#1}}

\newcommand\Dual{\mathcal D}
\newcommand\Duality\Dual

\newcommand\spinccan{{\mathfrak k}}

\newcommand\HFK{HFK}

\newcommand\HFKa{\widehat\HFK}

\begin{document}
\title{Holomorphic disks and genus bounds}
\asciiauthors{Peter Ozsvath and Zoltan Szabo}
\author{Peter Ozsv\'ath\\Zolt{\'a}n Szab{\'o}}
\coverauthors{Peter Ozsv\noexpand\'ath\\Zolt{\noexpand\'a}n Szab{\noexpand\'o}}
\address{Department of
Mathematics, Columbia University\\New York, NY 10025, USA\quad
{\rm and}\\Institute for Advanced Study, Princeton, New Jersey 08540, USA
\\{\rm and}\\Department of
Mathematics, Princeton University\\Princeton, New Jersey 08544, USA}

\gtemail{\mailto{petero@math.columbia.edu}\qua{\rm and}\qua
\mailto{szabo@math.princeton.edu}}
 
\asciiaddress{PO: Department of
Mathematics, Columbia University\\New York, NY 10025, USA
and\\Institute for Advanced Study, Princeton, New Jersey 08540, USA\\
SZ: Department of
Mathematics, Princeton University\\Princeton, New Jersey 08544, USA}

\asciiemail{petero@math.columbia.edu, szabo@math.princeton.edu}


\begin{abstract}
We prove that, like the Seiberg--Witten monopole homology, the Heegaard
Floer homology for a three-manifold determines its Thurston norm. As a
consequence, we show that knot Floer homology detects the genus of a
knot. This leads to new proofs of certain results previously obtained
using Seiberg--Witten monopole Floer homology (in collaboration with
Kronheimer and Mrowka). It also leads to a purely Morse-theoretic
interpretation of the genus of a knot.  The method of proof shows that
the canonical element of Heegaard Floer homology associated to a
weakly symplectically fillable contact structure is non-trivial. In
particular, for certain three-manifolds, Heegaard Floer homology gives
obstructions to the existence of taut foliations.
\end{abstract}

\asciiabstract{We prove that, like the Seiberg-Witten monopole
homology, the Heegaard Floer homology for a three-manifold determines
its Thurston norm. As a consequence, we show that knot Floer homology
detects the genus of a knot. This leads to new proofs of certain
results previously obtained using Seiberg-Witten monopole Floer
homology (in collaboration with Kronheimer and Mrowka). It also leads
to a purely Morse-theoretic interpretation of the genus of a knot.
The method of proof shows that the canonical element of Heegaard Floer
homology associated to a weakly symplectically fillable contact
structure is non-trivial. In particular, for certain three-manifolds,
Heegaard Floer homology gives obstructions to the existence of taut
foliations.}

\maketitle

\section{Introduction}

The purpose of this paper is to verify that the Heegaard Floer
homology of ~\cite{HolDisk} determines the Thurston semi-norm of its
underlying three-manifold.  This further underlines the relationship
between Heegaard Floer homology and Seiberg--Witten monopole Floer
homology of~\cite{KMbook}, 
for which an analogous result has been established by
Kronheimer and Mrowka, cf.~\cite{KMThurston}.

Recall that Heegaard Floer homology $\uHFa(Y)$ is a
finitely generated, $\Zmod{2}$--graded $\Z[H^1(Y;\Z)]$--module associated
to a closed, oriented three-manifold $Y$.  This group in turn admits a
natural splitting indexed by $\SpinC$ structures $\spinc$ over $Y$,
$$\uHFa(Y)=\bigoplus_{\spinc\in\SpinC(Y)}\uHFa(Y,\spinc).$$ (We
adopt here notation from~\cite{HolDisk}; the hat signifies here the
simplest variant of Heegaard Floer homology, while the underline
signifies that we are using the construction with ``twisted coefficients'',
cf.\ Section~8 of~\cite{HolDiskTwo}.)

The {\em Thurston semi-norm}~\cite{Thurston} on the two-dimensional
homology of $Y$ is the function $$\Theta\co H_2(Y;\Z)
\longrightarrow \Z^{\geq 0}$$ defined as follows.  The {\em
complexity} of a compact, oriented two-manifold $\chi_+(\Sigma)$ is
the sum over all the connected components $\Sigma_i\subset \Sigma$
with positive genus $g(\Sigma_i)$ of the quantity
$2g(\Sigma_i)-2$. The Thurston semi-norm of a homology class $\xi\in
H_2(Y;\Z)$ is the minimum complexity of any embedded representative of
$\xi$. (Thurston extends this function by linearity to a semi-norm
$\Theta\co H_2(Y;\Q) \longrightarrow \Q$.)

Our result now is the following:

\begin{theorem}
\label{thm:ThurstonNorm}
The $\SpinC$ structures $\spinc$ over $Y$ for which 
the Heegaard Floer homology $\uHFa(Y,\spinc)$ is
non-trivial determine the
Thurston semi-norm on $Y$, in the sense that:
$$\Theta(\xi)=
\max_{\{\spinc\in\SpinC(Y)\big| \uHFa(Y,\spinc)\neq 0\}} |\langle
c_1(\spinc), \xi\rangle |$$ for any $\xi\in H_2(Y;\Z)$.
\end{theorem}

The above theorem has a consequence for the ``knot Floer homology''
of~\cite{HolDiskKnots}, \cite{RasmussenThesis}.  For simplicity, we
state this for the case of knots in $S^3$.

Recall that knot Floer homology is a bigraded Abelian group
associated to an oriented knot $K\subset S^3$,
$$\HFKa(K)=\bigoplus_{d\in\Z, s\in \Z} \HFKa_d(K,s).$$ These groups
are a refinement of the Alexander polynomial of $K$, in the sense that
$$\sum_s \chi\left(\HFKa_*(K,s)\right) T^s=\Delta_K(T),$$ where here
$T$ is a formal variable, $\Delta_K(T)$ denotes the symmetrized
Alexander polynomial of $K$, and
$$\chi\left(\HFKa_*(K,s)\right)=\sum_{d\in\Z}(-1)^d
\Rk~\HFKa_d(K,s),$$ 
(cf.\ Equation~1 of~\cite{HolDiskKnots}).
One consequence of 
the proof of Theorem~\ref{thm:ThurstonNorm}
is the following quantitative sense in which $\HFKa$ distinguishes the
unknot:

\begin{theorem}
\label{thm:Knots}
Let $K\subset S^3$ be a knot, then the Seifert genus of $K$ is the
largest integer $s$ for which the group $\HFKa_*(K,s)\neq 0$.
\end{theorem}

This result in turn leads to an alternate proof of a theorem 
proved jointly by Kronheimer, Mrowka, and us~\cite{KMOSz}, first
conjectured by Gordon~\cite{GordonConjecture} (the cases where $p=0$
and $\pm 1$ follow from theorems of Gabai~\cite{GabaiKnots} and Gordon and
Luecke~\cite{GorLueckI} respectively):

\begin{cor}\label{cor:KMOSz}
{\rm\cite{KMOSz}}\qua
Let $K\subset S^3$ be a knot with the property that for some integer
$p$, $S^3_p(K)$ is diffeomorphic to $S^3_p(U)$ (where here $U$ is the
unknot) under an orientation-preserving diffeomorphism, then $K$ is
the unknot.
\end{cor}

The first ingredient in the proof of Theorem~\ref{thm:ThurstonNorm} is
a theorem of Gabai~\cite{Gabai} which expresses the minimal genus
problem in terms of taut foliations.  This result, together with a
theorem of Eliashberg and Thurston~\cite{EliashbergThurston} gives a
reformulation in terms of certain symplectically semi-fillable contact
structures.  The final breakthrough which makes this paper possible is
an embedding theorem of Eliashberg~\cite{Eliashberg}, see
also~\cite{Etnyre} and~\cite{OhtaOno}, which shows that a symplectic
semi-filling of a three-manifold can be embedded in a closed,
symplectic four-manifold. From this, we then appeal to a
theorem~\cite{HolDiskSymp}, which implies the non-vanishing of the
Heegaard Floer homology of a three-manifold which separates a closed,
symplectic four-manifold. This result, in turn, rests on the
topological quantum field-theoretic properties of Heegaard Floer
homology, together with the suitable handle-decomposition of an
arbitrary symplectic four-manifold induced from the Lefschetz pencils
provided by Donaldson~\cite{DonaldsonLefschetz}.  (The non-vanishing
result from~\cite{HolDiskSymp} is analogous to a non-vanishing theorem
for the Seiberg--Witten invariants of symplectic manifolds proved by
Taubes, cf.~\cite{TaubesSympI} and~\cite{TaubesSympII}.)

\subsection{Contact structures}

In another direction, the strategy of proof for
Theorem~\ref{thm:ThurstonNorm} shows that, just like its
gauge-theoretic counterpart, the Seiberg--Witten monopole Floer homology,
Heegaard Floer homology provides obstructions to the
existence of weakly symplectically fillable contact structures on a
given three-manifold, compare~\cite{KMcontact}.

For simplicity, we restrict attention now to the case where $Y$ is a
rational homology three-sphere, and hence $\HFa(Y)\cong
\uHFa(Y)$. In~\cite{Contact}, we constructed an invariant
$c(\xi)\in\HFa(Y)$, which we showed to be non-trivial for Stein
fillable contact structures. In Section~\ref{sec:HFContact}, we
generalize this to the case of symplectically semi-fillable contact
structures (see Theorem~\ref{thm:HFaWeakFilling} for a precise
statement).  It is very interesting to see if this non-vanishing result
can be generalized to the case of tight contact structures.
(Of course, in the
case where $b_1(Y)>0$, a reasonable formulation of this
question requires the use of twisted coefficients, cf.\
Section~\ref{sec:HFContact} below.)

In Section~\ref{sec:HFContact} 
we also prove a non-vanishing
theorem using the ``reduced Heegaard Floer homology'' $\HFpRed(Y)$
(for the image of $c(\xi)$ under a natural map $\HFa(Y) \longrightarrow
\HFpRed(Y)$), in the case where $b_2^+(W)>0$ or $W$ is a weak
symplectic semi-filling with more than one boundary
component. According to a result of Eliashberg and
Thurston~\cite{EliashbergThurston}, a taut foliation ${\mathcal F}$ on
$Y$ induces such a structure.

One consequence of this is an obstruction to the existence of such a
filling (or taut foliation) for a certain class of three-manifolds
$Y$. An $L$--space~\cite{NoteLens} is a rational homology three-sphere
with the property that $\HFa(Y)$ is a free $\Z$--module whose rank
coincides with the number of elements in $H_1(Y;\Z)$. Examples include
all lens spaces, and indeed all Seifert fibered spaces with positive
scalar curvature. More interesting examples are constructed as
follows: if $K\subset S^3$ is a knot for which $S^3_p(K)$ is an
$L$--space for some $p>0$, then so is $S^3_r(K)$ for all rational
$r>p$. A number of $L$--spaces are constructed in~\cite{NoteLens}. It
is interesting to note the following theorem of N{\'e}methi: a
three-manifold $Y$ is an $L$--space which is obtained as a plumbing of
spheres if and only if it is the link of a rational surface
singularity~\cite{Nemethi}. $L$--spaces in the context of
Seiberg--Witten monopole Floer homology are constructed in
Section~( of~\cite{KMOSz} (though the
constructions there apply equally well in the context of Heegaard
Floer homology).

The following theorem should be compared with~\cite{Lisca}, \cite{OhtaOno}
and~\cite{KMOSz} (see also~\cite{LiscaStipsicz}):

\begin{theorem}
\label{thm:Lspaces}
An $L$--space $Y$ has no symplectic semi-filling with disconnected boundary; 
and all its symplectic
fillings have $b_2^+(W)=0$. In particular, $Y$ admits no taut foliation.
\end{theorem}

\subsection{Morse theory and minimal genus}

Theorem~\ref{thm:ThurstonNorm} admits a reformulation which relates
the minimal genus problem directly in terms of Morse theory on the
underlying three-manifold.  For simplicity, we state this in the case
where $M$ is the complement of a knot $K\subset S^3$.

Fix a knot $K\subset S^3$. A perfect Morse function is said to be {\em
compatible with $K$}, if $K$ is realized as a union of two of the
flows which connect the index three and zero critical points (for some
choice of generic Riemannian metric $\mu$ on $S^3$).  Thus, the knot $K$ is
specified by a Heegaard diagram for $S^3$, equipped with two
distinguished points $w$ and $z$ where the knot $K$ meets the Heegaard
surface. In this case, a {\em simultaneous trajectory} is a collection
$\x$ of gradient flowlines for the Morse function which connect all
the remaining (index two and one) critical points of $f$.
From the point of view of Heegaard diagrams, a simultaneous
trajectory is an intersection point in the $g$--fold
symmetric product of $\Sigma$, $\Sym^g(\Sigma)$,
(where $g$ is the genus of $\Sigma$) of
two $g$--dimensional tori $\Ta=\alpha_1\times ...\times \alpha_g$ and
$\Tb=\beta_1\times...\times\beta_g$, where here $\{\alpha_i\}_{i=1}^g$
resp. $\{\beta_i\}_{i=1}^g$ denote the attaching circles of the two
handlebodies.

Let $X=X(f,\mu)$ denote the set of simultaneous trajectories. Any two
simultaneous trajectories differ by a one-cycle in the knot complement
$M$ and hence, if we fix an identification $H_1(M;\Z)\cong \Z$, we
obtain a difference map $$\epsilon \co X \times X \longrightarrow
\Z.$$ There is a unique map $s\co X \longrightarrow \Z$ with the
properties that $s(\x)-s(\y)=\epsilon(\x,\y)$ for all ${\mathbf x},
{\mathbf y}\in X$, and also $\#\{\x\big|s(\x)=i\}\equiv
\#\{\x\big|s(\x)=-i\}\pmod{2}$ for all $i\in\Z$.

Although we will not need this here, it is worth pointing out that
simultaneous trajectories can be viewed as a generalization of some
very familiar objects from knot theory. To this end, note that a knot
projection, together with a distinguished edge, induces in a natural
way a compatible Heegaard diagram. The simultaneous trajectories for
this Heegaard diagram can be identified with the ``Kauffman states''
for the knot projection; see~\cite{Kauffman} for an account of
Kauffman states, and~\cite{AltKnots} for their relationship with
simultaneous trajectories.

The following is a corollary of Theorem~\ref{thm:ThurstonNorm}.

\begin{cor}
\label{cor:MorseInterp}
The Seifert genus of a knot $K$ is the minimum over all compatible
Heegaard diagrams for $K$ of the maximum of $s(\x)$ over all the
simultaneous trajectories.
\end{cor}

It is very interesting to compare the above purely Morse-theoretic
characterization of the Seifert genus with Kronheimer and Mrowka's
purely differential-geometric characterization of the Thurston semi-norm on
homology in terms of scalar curvature, arising from the Seiberg--Witten
equations, cf.~\cite{KMThurston}. It would also be interesting to find
a more elementary proof of the above result.

\subsection{Remark}
This paper completely avoids the machinery of gauge theory and the
Seiberg--Witten equations. However, much of the general strategy
adopted here is based on the proofs of analogous results in monopole
Floer homology which were obtained by Kronheimer and Mrowka,
cf.~\cite{KMThurston}. It is also worth pointing out that although
the construction of Heegaard Floer homology is completely different
from the construction of Seiberg--Witten monopole Floer homology, the
invariants are conjectured to be isomorphic. (This conjecture should be
viewed in the light of the celebrated theorem of Taubes relating the
Seiberg--Witten invariants of closed symplectic manifolds with their
Gromov--Witten invariants, cf.~\cite{TaubesSWGr}.)

\subsection{Organization}
We include some preliminaries on contact geometry in
Section~\ref{sec:ContactPreliminaries}, and a quick review of Heegaard
Floer homology in Section~\ref{sec:Review}. In
Section~\ref{sec:HFContact}, we prove the non-vanishing results for
symplectically semi-fillable contact structures (including
Theorem~\ref{thm:Lspaces}). In Section~\ref{sec:ThurstonNorm} we turn
to the proofs of Theorems~\ref{thm:ThurstonNorm} and \ref{thm:Knots}
and 
Corollaries~\ref{cor:KMOSz} and \ref{cor:MorseInterp}.

\subsection{Acknowledgements}
This paper would not have been possible without the fundamental new
result of Yakov Eliashberg~\cite{Eliashberg}.  We would like to thank
Yasha for explaining his result to us, and for several illuminating
discussions. We would also like to thank Peter Kronheimer, Paolo
Lisca, Tomasz Mrowka, and Andr{\'a}s Stipsicz for many fruitful
discussions. We would especially like to thank Kronheimer and Mrowka
whose work in Seiberg--Witten monopole homology has served as an
inspiration for this paper.

PSO was partially supported by NSF grant numbers DMS-0234311,
DMS-0111298, and FRG-0244663.  ZSz was partially supported by NSF
grant numbers DMS-0107792 and FRG-0244663, and a Packard Fellowship.

\section{Contact geometric preliminaries}
\label{sec:ContactPreliminaries}

The three-manifolds we consider in this paper will always be oriented
and connected (unless specified otherwise).  A contact structure $\xi$
is a nowhere integrable two-plane distribution in $TY$. The contact
structures we consider in this paper will always be cooriented, and
hence (since our three-manifolds are also oriented) the two-plane
distributions $\xi$ are also oriented.  Indeed, they can be described
as the kernel of some smooth one-form $\alpha$ with the property that
$\alpha\wedge d\alpha$ is a volume form for $Y$ (with respect to its
given orientation). The form $d\alpha$ induces the orientation on
$\xi$.

A contact structure $\xi$ over $Y$ naturally gives rise to a $\SpinC$
structure, its {\em canonical $\SpinC$ structure}, written
$\spinccan(\xi)$, cf.~\cite{KMcontact}. Indeed, $\SpinC$ structures in
dimension three can be viewed as equivalence classes of nowhere
vanishing vector fields over $Y$, where two vector fields are
considered equivalent  if they are
homotopic in the complement of a ball in $Y$, cf.~\cite{Turaev},
\cite{GompfStipsicz}. Dually, an oriented two-plane distribution
gives rise to an equivalence class of nowhere vanishing vector fields
(which are transverse to the distribution, and form a positive basis
for $TY$).  Now, the canonical $\SpinC$ structure of a contact
structure is the $\SpinC$ structure associated to its two-plane
distribution.  The first Chern class of the canonical $\SpinC$
structure $\spinccan(\xi)$ is the first Chern class of $\xi$,
thought of now as a complex line bundle over $Y$.

Four-manifolds considered in this paper are also oriented.  A
symplectic four-manifold $(W,\omega)$ is a smooth four-manifold equipped with
a smooth two-form $\omega$ satisfying $d\omega=0$ and also the
non-degeneracy condition that $\omega\wedge \omega$ is a volume form
for $W$ (compatible with its given orientation).

Let $(W,\omega)$ be a compact, symplectic four-manifold $W$ with
boundary $Y$.  A four-manifold $W$ is said to have {\em convex
boundary} if there is a contact structure $\xi$ over $Y$ with the
property that the restriction of $\omega$ to the two-planes of $\xi$
is everywhere positive, cf.~\cite{EliashbergGromov}.  Indeed, if we
fix the contact structure $Y$ over $\xi$, we say that $W$ is a {\em
convex weak symplectic filling of $(Y,\xi)$}.  If $W$ is a convex weak
symplectic filling of a possibly disconnected three-manifold $Y'$ with contact
structure $\xi'$, and if $Y\subset Y'$ is a connected
subset with induced contact structure $\xi$, then we say that $W$ is a
{\em convex, weak semi-filling of $(Y,\xi)$}.  Of course, if a
symplectic four-manifold $W$ has boundary $Y$, equipped with a contact
structure $\xi$ for which the restriction of $\omega$ is everywhere
negative, we say that $W$ has {\em concave boundary}, and that $W$ is
a {\em concave weak symplectic filling of $Y$}.  (We use the term
``weak'' here to be consistent with the accepted terminology from
contact geometry. We will, however, never use the notion of strong
symplectic fillings in this paper.)

If a contact structure $(Y,\xi)$ admits a weak convex symplectic
filling, it is called {\em weakly fillable}. Note that every contact
structure $(Y,\xi)$ can be realized as the concave boundary of some
symplectic four-manifold (cf.~\cite{EtnyreHonda}, \cite{Gay},
and~\cite{Eliashberg}).  This is one justification for dropping the
modifier ``convex'' from the terminology ``weakly fillable''.  If a
contact structure $(Y,\xi)$ admits a weak symplectic semi-filling,
then it is called {\em weakly semi-fillable}. According to a recent
result of Eliashberg (cf.~\cite{Eliashberg}, restated in
Theorem~\ref{thm:Eliashberg} below) any weakly semi-fillable contact
structure is weakly fillable, as well. 

A symplectic structure $(W,\omega)$ endows $W$ with a canonical
$\SpinC$ structure, denoted $\spinccan(\omega)$,
cf.~\cite{TaubesSympI}. This can be thought of as the canonical
$\SpinC$ structure associated to any almost-complex structure $J$ over
$W$ compatible with $\omega$, compare~\cite{TaubesSympI}. In
particular, the first Chern class the $\SpinC$ structure
$\spinccan(\omega)$ is the first Chern class of its complexified
tangent bundle.  If $(W,\omega)$ has convex boundary $(Y,\xi)$, then
the restriction of the canonical $\SpinC$ structure over $W$ to $Y$ is
the canonical $\SpinC$ structure of the contact structure $\xi$.

\subsection{Foliations and contact structures}

Recall that a taut foliation is a foliation ${\mathcal F}$ which comes
with a two-form $\omega$ which is positive on the leaves of ${\mathcal
F}$ (note that like our contact structures, all the foliations we
consider here are cooriented and hence oriented). An {\em irreducible
three-manifold} is a three-manifold $Y$ with $\pi_2(Y)=0$. A
fundamental result of Gabai states that if $Y$ is irreducible and
$\Sigma_0\subset Y$ is an embedded surface which minimizes complexity
in its homology class, and with has no spherical or toroidal
components, then there is a smooth, taut foliation ${\mathcal F}$
which contains $\Sigma_0$ as a union of compact leaves. In particular,
this shows that if $Y$ is an irreducible three-manifold with
non-trivial Thurston semi-norm, and $\Sigma\subset Y$ is an embedded
surface which minimizes complexity in its homology class, then there
is a smooth, taut foliation ${\mathcal F}$ with the property that $\langle
c_1({\mathcal F}),[\Sigma]\rangle = -\chi_+(\Sigma)$.  (Here, we let
${\mathcal F}$ be a taut foliation whose closed leaves include all the
components of $\Sigma$ with genus greater than one.)

The link between taut foliations and semi-fillable contact structures
is provided by an observation of Eliashberg and Thurston,
cf.~\cite{EliashbergThurston}, according to which if $Y$ admits a
smooth, taut foliation ${\mathcal F}$, then $W=[-1,1]\times Y$ can be
given the structure of a convex symplectic manifold, where here the
two-plane fields $\xi_\pm$ over $\{\pm 1\}\times Y$ are homotopic to
the two-plane field of tangencies to ${\mathcal F}$.

\section{Heegaard Floer homology}
\label{sec:Review}

Heegaard Floer homology is a collection of $\Zmod{2}$--graded homology
theories associated to three-manifolds, which are functorial under
smooth four-dimensional cobordisms (cf.~\cite{HolDisk} for their
constructions, and~\cite{HolDiskFour} for the verification of their
functorial properties).

There are four variants, $\HFa(Y)$, $\HFm(Y)$, $\HFinf(Y)$, and
$\HFp(Y)$.  $\HFm(Y)$ is the homology of a complex over the polynomial
ring $\Z[U]$, $\HFinf(Y)$ is the associated ``localization'' (i.e.\ it
is the homology of the complex associated to tensoring with the ring
of Laurent polynomials over $U$), $\HFp(Y)$ is associated to the
cokernel of the localization map, and finally $\HFa(Y)$ is the homology of
the complex
associated to setting $U=0$.  Indeed, all these groups admit splittings
indexed by $\SpinC$ structures over $Y$. The various groups are 
related by long exact sequences
\begin{equation}
\label{eq:ExactSequences}
\begin{CD}
...@>>>\HFa(Y,\spinct)@>{i}>>\HFp(Y,\spinct)@>{U}>>\HFp(Y,\spinct)@>>> ...\\
...@>>>\HFm(Y,\spinct)@>{j}>>\HFinf(Y,\spinct)@>{\pi}>>\HFp(Y,\spinct)@>>>...,
\end{CD}
\end{equation} 
where here $\spinct\in\SpinC(Y)$.  The ``reduced Heegaard Floer
homology'' $\HFpRed(Y,\spinct)$ is the cokernel of the map
$\pi$. Sometimes we distinguish this from $\HFmRed(Y,\spinct)$, which
is the kernel of the map $j$, though these two $\Z[U]$ modules are
identified in the long exact sequence above.

For $Y=S^3$, we have that $\HFa(S^3)\cong \Z$. We can now lift the
$\Zmod{2}$ grading to an absolute $\Z$--grading on all the groups,
using the following conventions. The group $\HFa(S^3)\cong \Z$ is
supported in dimension zero, the maps $i$, $j$, and $\pi$ from
Equation~\eqref{eq:ExactSequences} preserve degree, and $U$ decreases
degree by two. Indeed, for $S^3$, we have an identification of $\Z[U]$
modules: $$
\begin{CD}
0@>>>\HFm(S^3) @>>>\HFinf(S^3)@>>>\HFp(S^3)@>>>0 \\
&&@V{=}VV @V{=}VV @V{=}VV \\
0@>>>U\cm \Z[U] @>>> \Z[U,U^{-1}] @>>>\Z[U,U^{-1}]/U\cm\Z[U]@>>>0,
\end{CD}$$
where here the element $1\in\Z[U,U^{-1}]$ lies in grading zero and $U$ decreases
grading by two. (See~\cite{AbsGraded} for a definition of absolute gradings
in more general settings.)

To state functoriality, we must first discuss maps associated to
cobordisms. Let $W_1$ be a smooth, oriented four-manifold with
$\partial W_1 = -Y_1 \cup Y_2$, where here $Y_1$ and $Y_2$ are
connected. (Here, of course, $-Y_1$ denotes the three-manifold
underlying $Y_1$, endowed with the opposite orientation.) In this
case, we sometimes write $W_1\co Y_1\longrightarrow Y_2$; or,
turning this around, we can view the same four-manifold as giving a
cobordism $W_1\co -Y_2\longrightarrow -Y_1$. There is an associated
map $$\Fa{W_1}\co \HFa(Y_1)
\longrightarrow \HFa(Y_2),$$ 
well-defined up to an overall multiplication by $\pm 1$,
which can be decomposed along $\SpinC$
structures over $W_1$: $$\Fa{W_1,\spinc}\co
\HFa(Y_1,\spinct_1) \longrightarrow \HFa(Y_2,\spinct_2),$$ where here
$\spinct_i=\spinc|_{Y_i}$, i.e.\ so that
$$\Fa{W_1}=\sum_{\spinc\in\SpinC(W_1)}\Fa{W_1,\spinc}.$$ There are similarly
induced maps $\Fp{W_1,\spinc}$ on $\HFp$ which are equivariant under the
action of $\Z[U]$. For $\HFinf$ and $\HFm$, there are again induced
maps $\Finf{W_1,\spinc}$ and $\Fm{W_1,\spinc}$ for each fixed $\SpinC$
structure $\spinc\in \SpinC(W_1)$ (but now, we can no longer sum maps
over all $\SpinC$ structures, since infinitely many might be non-trivial). 
Indeed, these maps are compatible
with the natural maps from Diagram~\eqref{eq:ExactSequences}; for example,
all the squares in the following diagram commute:
$$\begin{CD}
...@>>>\HFm(Y_1,\spinct_1)@>>>\HFinf(Y_1,\spinct_1)@>>>\HFp(Y_1,\spinct_1)@>>>...\\
&& @V{\Fm{W_1,\spinc}}VV @V{\Finf{W_1,\spinc}}VV @V{\Fp{W_1,\spinc}}VV \\
...@>>>\HFm(Y_2,\spinct_2)@>>>\HFinf(Y_2,\spinct_2)@>>>\HFp(Y_2,\spinct_2)@>>>...
\end{CD}$$

Functoriality of Floer homology is to be interpreted in the following
sense. Let $W_1\co Y_1\longrightarrow Y_2$ and $W_2 \co
Y_2\longrightarrow Y_3$. We can form then the composite cobordism
$$W_1\#_{Y_2}W_2\co Y_1\longrightarrow Y_3.$$ We claim that
for each $\spinc_i\in\SpinC(W_i)$ with $\spinc_1|_{Y_{2}}=\spinc_2|_{Y_2}$,
we have that
\begin{equation}
\label{eq:ProductFormula}
\sum_{\{\spinc\in\SpinC(W_1\#_{Y_2}W_2) \big|\spinc|_{W_i}=\spinc_i\}}
\Fa{W,\spinc} = \Fa{W_2,\spinc_2}\circ \Fa{W_1,\spinc_1},
\end{equation}
with analogous formulas for $\HFm$, $\HFinf$, and $\HFp$ as well
(this is
the ``composition law'',
Theorem~3.4 of~\cite{HolDiskFour}).

Of these theories, $\HFinf$ is the weakest at distinguishing manifolds.
For example, if
$W\co Y_1\longrightarrow Y_2$ is a cobordism with $b_2^+(W)>0$,
then for any $\SpinC$ structure $\spinc\in\SpinC(W)$
the induced map
$$\Finf{W,\spinc}\co \HFinf(Y_1,\spinc|_{Y_1}) \longrightarrow 
\HFinf(Y_2,\spinc|_{Y_2})$$
vanishes (cf.\ Lemma~8.2
of~\cite{HolDiskFour}).

Floer homology can be used to construct an invariant for smooth
four-manifolds $X$ with $b_2^+(X)>1$ (here, $b_2^+(X)$ denotes the
dimension of the maximal subspace of $H^2(X;\R)$ on which the
cup-product pairing is positive-definite) endowed with a $\SpinC$
structure $\spinc\in\SpinC(X)$ $$\Phi_{X,\spinc}\co
\Z[U]\longrightarrow \Z,$$ which is well-defined up to an overall
sign. This invariant is analogous to the Seiberg--Witten invariant,
cf.~\cite{Witten}. This map is a homogeneous element in $\Hom(\Z[U],\Z)$
with degree given by
$$\frac{c_1(\spinc)^2-2\chi(X)-3\sigma(X)}{4}.$$ For a
fixed four-manifold $X$, the invariant $\Phi_{X,\spinc}$ is
non-trivial for only finitely many $\spinc\in\SpinC(X)$.  (Note that
the four-manifold invariant $\Phi_{X,\spinc}$ constructed in~\cite{HolDiskFour}
is slightly more general, as it incorporates the action of
$H_1(X;\Z)$, but we do not need this extra structure for our present
applications.)

The invariant is constructed as follows. Let $X$ be a four-manifold, and
fix a separating hypersurface $N\subset X$ with $0=\delta
H^1(N;\Z)\subset H^2(X;\Z)$, so that $X=X_1 \cup_N X_2$, with
$b_2^+(X_i)>0$ for $i=1,2$. (Here, $\delta\co H^1(Y;\Z)\longrightarrow H^2(X;\Z)$
is the connecting homomorphism in the Mayer-Vietoris sequence
for the decomposition of $X$ into $X_1$ and $X_2$.)
Such a separating three-manifold is called an {\em
admissible cut} in the terminology of~\cite{HolDiskFour}. Given such a cut,
delete balls $B_1$ and $B_2$ from 
$X_1$ and $X_2$ respectively, and consider the diagram: $$
\begin{CD}
&&&&\HFm(S^3)@>>>\HFinf(S^3) \\
&&&&@V{\Fm{X_1-B_1,\spinc_1}}VV  @V{\Finf{X_1-B_1,\spinc_1}}V{0}V \\
\HFinf(N,\spinct) @>>>\HFp(N,\spinct)@>>>\HFm(N,\spinct)@>>>\HFinf(N,\spinct) \\
@V{0}V{\Finf{X_2-B_2,\spinc_2}}V @VV{\Fp{X_2-B_2,\spinc_2}}V \\
\HFinf(S^3)@>>>\HFp(S^3), 
\end{CD}$$
where here $\spinct=\spinc|_{N}$ and $\spinc_i=\spinc|_{X_i}$.  Since
the two maps indicated with $0$ vanish (as $b_2^+(X_i-B_i)>0$), there
is a well-defined map $$\Fmix{X-B_1-B_2,\spinc}\co\HFm(S^3)\longrightarrow
\HFp(S^3),$$ which factors through $\HFpRed(N,\spinct)$. 

The invariant $\Phi_{X,\spinc}$ corresponds to
$\Fmix{X-B_1-B_2,\spinc}$ under the natural identification
$$\Hom_{\Z[U]}(\Z[U],\Z[U,U^{-1}]/\Z[U])\cong \Hom(\Z[U],\Z)$$
According to Theorem~9.1
of~\cite{HolDiskFour}, $\Phi_{X,\spinc}$ is a smooth four-manifold
invariant.

The following property of the invariant is immediate from its
definition: if $X=X_1\cup_N X_2$ where $N$ is a rational homology
three-sphere with $\HFpRed(N)=0$, and the four-manifolds $X_i$ have
the property that $b_2^+(X_i)>0$, then for each $\spinc\in\SpinC(X)$,
$$\Phi_{X,\spinc}\equiv 0.$$ 
The second property which we rely on heavily in this paper is the
following analogue of a theorem of Taubes~\cite{TaubesSympI} and
\cite{TaubesSympII} for the Seiberg--Witten invariants for
four-manifolds: if $(X,\omega)$ is a smooth, closed, symplectic
four-manifold with $b_2^+(X)>1$, then if
$\spinccan(\omega)\in\SpinC(X)$ denotes its canonical $\SpinC$
structure, then we have that $$\Phi_{X,\spinccan(\omega)}\equiv \pm
1,$$ while if $\spinc\in\SpinC(X)$ is any $\SpinC$ structure for which
$\Phi_{X,\spinc}\not\equiv 0$, then we have that 
$$\langle c_1(\spinccan(\omega))\cup
\omega,[X]\rangle \leq
\langle c_1(\spinc)\cup
\omega,[X]\rangle,$$
with equality iff $\spinc=\spinccan(\omega)$.
This result is 
Theorem~1.1
of~\cite{HolDiskSymp}, and its proof relies on a
combination of techniques from Heegaard Floer homology (specifically,
the surgery long exact sequence from~\cite{HolDiskTwo}) and
Donaldson's Lefschetz pencils for symplectic
manifolds,~\cite{DonaldsonLefschetz}.

\subsection{Three-manifolds with $b_1(Y)>0$}

There is a version of Floer homology with ``twisted coefficients''
which is relevant in the case where $b_1(Y)>0$. Fundamental to this
construction is a chain complex $\uCFa(Y)$ (and also corresponding
complexes $\uCFm$, $\uCFinf$, and $\uCFp$) with coefficients
in $\Z[H^1(Y;\Z)]$ which is a lift of the complex $\CFa(Y)$ (whose
homology calculates $\HFa(Y)$), in the following sense. Let $\Z$ be
the module over $\Z[H^1(Y;\Z)]$, where the elements of $H^1(Y;\Z)$ act
trivially. Then, there is an identification $\CFa(Y)\cong
\uCFa(Y)\otimes_{\Z[H^1(Y;\Z)]} \Z$. Thus, there is a change
of coefficient spectral sequences which relates the homology of
$\uCFa(Y)$, written $\uHFa(Y)$, with $\HFa(Y)$.

Indeed, given any module $M$ over $\Z[H^1(Y;\Z)]$, we can form the
group $$\uHFa(Y;M)=H_*\left(\uCFa(Y)\otimes_{\Z[H^1(Y;\Z)]}
M\right),$$ which gives Floer homology with coefficients twisted by
$M$.  The analogous construction in the other versions of Floer
homology gives groups $\uHFm(Y;M)$, $\uHFinf(Y;M)$, and
$\uHFp(Y;M)$. All of these are related by exact sequences analogous to
those in Diagram~\eqref{eq:ExactSequences}. In particular, we can
form a reduced group $\uHFpRed(Y;M)$, which is the cokernel of
the localization map $\uHFinf(Y;M)\longrightarrow \uHFp(Y;M)$.

In particular, if we fix a two-dimensional cohomology class
$[\omega]\in H^2(Y;\R)$, we can view $\Z[\R]$ as a module over
$\Z[H^1(Y;\Z)]$ via the ring homomorphism $$[\gamma]\mapsto T^{\int_Y
[\gamma]\wedge \omega}$$ (where here $T^r$ denotes the group-ring
element associated to the real number $r$). This gives us a notion of
twisted coefficients which we denote by $\uHFa(Y;[\omega])$.

This can be thought of explicitly as follows. Choose a Morse function
on $Y$ compatible with a Heegaard decomposition
$(\Sigma,\alphas,\betas,z)$, and fix also a two-cocycle $\omega$ over
$Y$ which represents $[\omega]$. We obtain a map from Whitney disks
$u$ in $\Sym^g(\Sigma)$ (for $\Ta$ and $\Tb$)
to two-chains in $Y$: $u$ induces a two-chain in $\Sigma$ with
boundaries along the $\alphas$ and $\betas$. These boundaries are
then coned off by following gradient trajectories for the $\alpha$--
and $\beta$--circles. Since $\omega$ is a cocycle, the evaluation of
$\omega$ on $u$ depends only on the homotopy class $\phi$ 
of $u$. We denote this evaluation by $\int_{[\phi]}
\omega$. (This determines an additive assignment in the terminology of
Section~8 of~\cite{HolDiskTwo}.)
The differential on $\uHFp(Y;[\omega])$
is given by
$$ {\underline\partial}^+ [\x,i] = \sum_{\y\in\Ta\cap\Tb} \sum_{\{\phi\in\pi_2(\x,\y)\big|\Mas(\phi)=1\}}
\#\left(\frac{\ModFlow(\phi)}{\R}\right)\cm 
T^{\int_{[\phi]}\omega} \cm [\y,i-n_z(\phi)],$$ where here we adopt
notation from~\cite{HolDiskTwo}: $\pi_2(\x,\y)$ denotes the space of
homotopy classes of Whitney disks in $\Sym^g(\Sigma)$ for $\Ta$ and
$\Tb$ connecting $\x$ and $\y$, $\Mas(\phi)$ denotes the formal
dimension of its space $\ModFlow(\phi)$ of holomorphic
representatives, and $n_z(\phi)$ denotes the intersection number of
$\phi$ with the subvariety $\{z\}\times \Sym^{g-1}(\Sigma)\subset
\Sym^g(\Sigma)$.

Now, if $W\co Y_1\longrightarrow Y_2$, and $M_1$ is a module 
over $H^1(Y_1;\Z)$, there is an induced map $$\uFp{W;M_1}\co \uHFp(Y_1,M_1)
\longrightarrow \uHFp(Y_2,M_1\otimes_{H^1(Y_1;\Z)} H^2(W, Y_1\cup
Y_2)),$$
well-defined up to the action by some unit in $\Z[H^2(Y_1\cup Y_2;\Z)]$,
defined as in
Subsection~3.1~\cite{HolDiskFour}.
(Indeed, in that discussion, the construction is separated according
to $\SpinC$ structures over $W$, which we drop at the moment for
notational simplicity.) In the case of $\omega$--twisted coefficients,
this gives rise to a map
$$\uFp{W;[\omega]}\co \uHFp(Y_1;[\omega]|_{Y_1})\longrightarrow
\uHFp(Y_2;[\omega]|_{Y_2})$$
(again, well-defined up to multiplication by $\pm T^c$ for some $c\in\R$)
which can be concretely described as follows.

Suppose for simplicity that $W$ is represented as a two-handle
addition, so that there is a corresponding ``Heegaard triple''
$(\Sigma,\alphas,\betas,\gammas,z)$. The corresponding four-manifold
$X_{\alpha,\beta,\gamma}$ represents $W$ minus a one-complex. Fix now
a two-cocycle $\omega$ representing $[\omega]\in H^2(W;\R)$. Again, a Whitney
triangle $u$ in $\Sym^g(\Sigma)$ for $\Ta$, $\Tb$, and $\Tc$ (with
vertices at $\x$, $\y$, and $\w$) determines a two-chain in
$X_{\alpha,\beta,\gamma}$, whose evaluation on $\omega$ depends on $u$
only through its induced homotopy class $\psi$ in $\pi_2(\x,\y,\w)$, denoted by
$\int_{[\psi]}\omega$. Now,
\begin{equation}
\label{eq:DefFp}
\uFp{W;[\omega]}[\x,i]
=\sum_{\y\in\Ta\cap \Tc}\sum_{\{\psi\in\pi_2(\x,\Theta,\y)\big|
\Mas(\psi)=0\}}\#\left(\ModFlow(\psi)\right)\cm T^{\int_{[\psi]}\omega}
\cm [\y,i-n_z(\psi)],
\end{equation}
where $\Theta\in\Tb\cap\Tc$ represents a canonical generator for the
Floer homology 
$\HF^{\leq}=H_*(U^{-1}\cdot {\mathrm CF}^-)$
of the three-manifold determined by
$(\Sigma,\betas,\gammas,z)$, which is a connected sum
$\#^{g-1}(S^2\times S^1)$. This can be extended to 
arbitrary (smooth, connected) cobordisms from $Y_1$ to $Y_2$ as
in~\cite{HolDiskFour}. 

(In the present discussion, since we have suppressed $\SpinC$
structures from the notation, a subtlety arises. The expression
analogous to Equation~\eqref{eq:DefFp}, only using $\HFm$, is not
well-defined since, in principle, there might be infinitely many
different homotopy classes which induce non-trivial maps -- i.e.\ we
are trying to sum the maps on $\HFm$ induced by infinitely many
different $\SpinC$ structures. However, if the cobordism $W$ has
$b_2^+(W)>0$, then there are only finitely many $\SpinC$ structures
which induce non-zero maps, according to
Theorem~3.3 of~\cite{HolDiskFour}.)

Note that when $W$ is a cobordism between two integral homology
three-spheres, the above construction is related to the construction
in the untwisted case by the formula
$$\uFp{W;[\omega]} = \pm T^{c}\cm \sum_{\spinc\in\SpinC(W)} 
T^{\langle c_1(\spinc)\cup [\omega], [W]\rangle} \cm \Fp{W,\spinc}$$
for some constant $c\in\R$.

\section{Invariants of weakly fillable contact structures}
\label{sec:HFContact}

We briefly review the construction here of the Heegaard Floer homology
element associated to a contact structure $\xi$ over the
three-manifold $Y$, $c(\xi)\in\HFa(-Y)$. After sketching the
construction, we describe  a refinement which lives in Floer
homology with twisted coefficients. 

The contact invariant is constructed with the help of some work of
Giroux.  Specifically, in~\cite{Giroux}, Giroux shows that contact
structures over $Y$ are in one-to-one correspondence with equivalence
classes of open book decompositions of $Y$, under an equivalence
relation given by a suitable notion of stabilization. Indeed, after
stabilizing, one can realize the open book with connected binding, and
with genus $g>1$ (both are convenient technical devices). In
particular, performing surgery on the binding, we obtain a cobordism
(obtained by a single two-handle addition) $W_0\co Y
\longrightarrow Y_0$, where here the three-manifold $Y_0$ fibers
over the circle. We call this cobordism a {\em Giroux two-handle}
subordinate to the contact structure over $Y$.  This cobordism is used
to construct $c(\xi)$, but to describe how, we must discuss the
Heegaard Floer homology for three-manifolds which fiber over the
circle.

Let $Z$ be a (closed, oriented) three-manifold endowed with the
structure of a fiber bundle $\pi\co Z\longrightarrow S^1$.  This
structure endows $Z$ with a canonical $\SpinC$ structure
$\spinccan(\pi)\in\SpinC(Z)$ (induced by the two-plane distribution
of tangents to the fiber of $\pi$). 
According to~\cite{HolDiskSymp}, if the genus
$g$ of the fiber is greater than one, then $$\HFp(Z,\spinccan(\pi))\cong
\Z.$$ In particular, there is a homogeneous generator $c_0(\pi)$ for
$\HFa(Z,\spinccan(\pi))\cong \Z\oplus \Z$ which maps to the generator 
$c^+_0(\pi)$ of
$\HFp(Z,\spinccan(\pi))$. This generator is, of course, uniquely determined
up to sign.

With these remarks in place, we can give the definition of the
invariant $c(\xi)$ associated to a contact structure over $Y$.  If $Y$
is given a contact structure, fix a compatible open book decomposition
(with connected binding, and fiber genus $g>1$), and consider the
corresponding Giroux two-handle $W_0\co -Y_0\longrightarrow -Y$
(which we have ``turned around'' here), and let $${\widehat
F}_{W_0}\co \HFa(-Y_0) \longrightarrow \HFa(-Y)$$ be the induced
map. Then, define $c(\xi)\in\HFa(-Y)/\{\pm 1\}$ to be the image
${\widehat F}_{W_0}(c_0(\pi))$. It is shown in~\cite{Contact} that
this element is uniquely associated (up to sign) to the contact
structure, i.e.\ it is independent of the choice of compatible open
book.  In fact, the element $c(\xi)$ is supported in the summand
$\HFa(Y,\spinccan(\xi))\subset\HFa(Y)$, where here $\spinccan(\xi)$ is
the canonical $\SpinC$ structure associated to the contact structure
$\xi$, in the sense described in
Section~\ref{sec:ContactPreliminaries}. (In particular, the canonical
$\SpinC$ structure of the fibration structure on $-Y_0$ is $\SpinC$
cobordant to the canonical $\SpinC$ structure of the contact structure
over $-Y$ via the Giroux two-handle.)

With the help of Giroux's characterization of Stein fillable contact
structures, it is shown in~\cite{Contact} that $c(\xi)$ is non-trivial
for a Stein structure. This non-vanishing result can be strengthened
considerably with the help of the following result of
Eliashberg~\cite{Eliashberg}.

\begin{theorem}[Eliashberg~\cite{Eliashberg}]
\label{thm:Eliashberg}
Let $(Y,\xi)$ be a contact three-manifold, which is the convex
boundary of some symplectic four-manifold $(W,\omega)$. Then, any
Giroux two-handle $W_0\co Y\longrightarrow Y_0$ can be completed to
give a compact symplectic manifold $(V,\omega)$ with concave boundary
$\partial (V,\omega)=(Y,\xi)$, so that $\omega$ extends smoothly over
$X=W\cup_Y V$.
\end{theorem}

Although Eliashberg's is the construction we need, concave fillings
have been constructed previously in a number of different contexts,
see for example~\cite{LiscaMatic}, \cite{AkbulutOzbagci},
\cite{EtnyreHonda}, \cite{Gay}, \cite{OhtaOno}. Indeed, since the
first posting of the present article, Etnyre pointed out to us an
alternate proof of Eliashberg's theorem~\cite{Etnyre},
see also \cite{OhtaOno}.

In the construction, $V$ is given as the union of the Giroux
two-handle with a surface bundle $V_0$ over a surface-with-boundary
which extends the fiber bundle structure over $Y_0$.  Moreover, the
fibers of $V_0$ are symplectic. By forming a symplectic sum if
necessary, one can arrange for $b_2^+(V)$ to be arbitrarily large.

To state the stronger non-vanishing theorem, we use a refinement of the contact
element using twisted coefficients. We can repeat the construction of
$c(\xi)$ with coefficients in any module $M$ over $\Z[H^1(Y;\Z)]$
(compare Remark~4.5 of~\cite{Contact}), to get
an element
$$c(\xi;M)\in\uHFa(Y;M)/\Z[H^1(Y;\Z)]^\times.$$ As the notation
suggests, this is an element $c(\xi;M)\in\uHFa(Y;M)$, which is
well-defined up to overall multiplication by a unit in the group-ring
$\Z[H^1(Y;\Z)]$.
Let $c^+(\xi;M)$ denote the image of $c(\xi;M)$ under the natural
map $\uHFa(-Y;M) \longrightarrow \uHFp(-Y;M)$, and let
$c^+_{\red}(\xi;M)$ denote its image under the projection\break
$\uHFp(-Y;M) \longrightarrow \uHFpRed(-Y;M)$.

In our applications, we will typically take the module $M$ to be
$\Z[\R]$, with the action specified by some two-form $\omega$ over
$Y$, so that we get $c(\xi;[\omega])\in\uHFa(-Y;[\omega])$.  The
following theorem should be compared with a theorem of Kronheimer and
Mrowka~\cite{KMcontact}, see also
Section~6 of~\cite{KMOSz}:

\begin{theorem}
\label{thm:HFaWeakFilling}
Let $(W,\omega)$ be a weak filling of a contact structure $(Y,\xi)$.
Then, the associated contact invariant $c(\xi;[\omega])$ is
non-trivial.  Indeed, it is non-torsion and primitive (as is its image
in $\uHFp(Y;[\omega])$.  Indeed, if $(W,\omega)$ is a
weak-semi-filling of $(Y,\xi)$ with disconnected boundary or
$(W,\omega)$ is a weak filling of $Y$ with $b_2^+(W)>0$, then the
reduced invariant $c^+_{\red}(\xi;[\omega])$ is non-trivial (and
indeed non-torsion and primitive).
\end{theorem}

\begin{proof}
Let $(W,\omega)$ be a symplectic filling of $(Y,\xi)$ with convex boundary.

Consider Eliashberg's cobordism bounding $Y$, $V=W_0\cup_{Y_0} V_0$,
where here $W_0\co Y \longrightarrow Y_0$ is the Giroux two-handle
and $V_0$ is a surface bundle over a surface-with-boundary. Now, the
union $$X=V_0\cup_{-Y_0}\cup W_0\cup_{-Y} W$$ is a closed, symplectic
four-manifold.  (As the notation suggests, we have ``turned around''
$W_0$, to think of it as a cobordism from $-Y_0$ to $-Y$; similarly
for $V_0$.) Arrange for $b_2^+(V_0)>1$, and decompose $V_0$ further by
introducing an admissible cut by $N$. Now, $N$ decompose $X$
into two pieces $X=X_1 \cup_N X_2$, where $b_2^+(X_i)>0$, and we can
suppose now that $X_2$ contains the Giroux cobordism, i.e.
\begin{equation}
\label{eq:DecomposeX2}
X_2 = (V_0-X_1) \cup_{-Y_0}\cup W_0 \cup_{-Y} W.
\end{equation}
Now, by the definition of $\Phi$, for any given $\spinc\in\SpinC(X)$,
there is an element $\theta\in \HFp(N,\spinc|_{N})$ 
with the property
that $$\Phi_{X,\spinc}=\Fp{X_2-B_2}(\theta).$$ 
(By definition of $\Phi$, the element $\theta$ here is any element of $\HFp(N,\spinc|_{N})$
whose image under the connecting homomorphism in
the second exact sequence in Equation~\eqref{eq:ExactSequences} coincides with 
the image of a generator of $\HFm(S^3)$ under the map $\Fm{X_1-B_1}\co \HFm(S^3) \longrightarrow \HFm(N,\spinc|_{N})$.)
Applying the product
formula for the decomposition of Equation~\eqref{eq:DecomposeX2}, we
get that $$\sum_{\eta\in H^1(Y;\Z)}
\Phi_{X,\spinccan(\omega)+\delta \eta}
= \Fp{W-B_2}\circ \Fp{W_0}\circ\Fp{V_0-X_1}(\theta).$$
In terms of $\omega$--twisted coefficients, we have that
$$\sum_{\eta\in H^1(Y_0;\Z)}\Phi_{X,\spinccan(\omega)+\delta \eta}\cm 
T^{\langle \omega\cup c_1(\spinccan(\omega)+\delta\eta),[X]\rangle}
= 
\uFp{W-B_2;[\omega]}\circ\uFp{W_0;[\omega]}\circ \uFp{V_0-X_1;[\omega]}(\underline\theta).$$
(Here, ${\underline\theta}\in\uHFp(N,\spinc|_{N};[\omega])$ is the analogue
of the class $\theta$ considered earlier.)
But $\HFp(Y_0,\spinct)\cong \Z[\R]$ is generated by $c^+_0(\pi)$
(where here $\pi\co Y_0\longrightarrow S^1$ is the projection
obtained from restricting the bundle structure over $V_0$, and 
$\spinct$ is the restriction of $\spinccan(\omega)$ to $Y_0$), 
so there is some element
$p(T)\in\Z[\R]$ with the property that
$\uFp{V_0-\nbd{F}}({\underline\theta})=p(T)\cm c^+(\pi)$. Thus, 
$$\sum_{\eta\in H^1(Y_0;\Z)}\Phi_{X,\spinccan(\omega)+\delta \eta}\cm 
T^{\langle \omega\cup c_1(\spinccan(\omega)+\delta\eta),[X]\rangle}=p(T)\cm 
\uFp{W-B_2}(c^+(\xi;[\omega])).$$
The left-hand-side here gives a polynomial in $T$ (well defined up to
an overall sign and multiple of $T$) whose lowest-order term is one,
according to Theorem~1.1
of~\cite{HolDiskSymp} (recalled in Section~\ref{sec:Review}). It
follows at once that $\uFp{W-B_2}(c^+(\xi;[\omega]))$ is
non-trivial. Indeed, it also follows that
$\uFp{W-B_2}(c^+(\xi;[\omega]))$ is a primitive homology class
(since the leading coefficient is $1$), and no multiple of
it zero. This implies the same for $c(\xi;[\omega])$.

Now, when $b_2^+(W)>0$, we use $Y$ as a cut for $X$ to show that the
induced element $c^+_{\red}(\xi;[\omega])$ is non-trivial (primitive
and torsion). In the case where $Y$ is semi-fillable with disconnected
boundary, we can close off the remaining boundary components as in
Theorem~\ref{thm:Eliashberg} to construct a new symplectic filling
$W'$ of $Y$ with one boundary component and $b_2^+(W')>0$, reducing to
the previous case.
\end{proof}

\medskip
\noindent{\bf{Proof of Theorem~\ref{thm:Lspaces}}}\qua A
three-manifold $Y$ is an $L$--space if it is a rational homology
three-sphere and $\HFa(Y)$ is a free $\Z$--module of rank
$|H_1(Y;\Z)|$. Note that for an $L$-space, $\HFpRed(Y)\otimes_\Z \Q =
0$.  This is an easy application of the long exact
sequence~\eqref{eq:ExactSequences}, together with the fact that the
the intersection of the kernel of $U\co \HFp(Y) \longrightarrow
\HFp(Y)$ with the image of $\HFinf(Y)$ inside $\HFp(Y)$ has rank
$|H_1(Y;\Z)|$, since $\HFinf(Y)\cong \Z[U,U^{-1}]$ (cf.\ Theorem~10.1
of~\cite{HolDiskTwo}), the map from $\HFinf(Y)$ to $\HFp(Y)$ is an
isomorphism in all sufficiently large degrees (i.e.\ $U^{-n}$ for $n$
sufficiently large), and it is trivial in all sufficiently small
degrees.

For a three-manifold $Y$ with $b_1(Y)=0$, $\uHFp(Y;[\omega])\cong
\HFp(Y)\otimes_{\Z}\Z[\R]$, since $[\omega]\in H^2(Y;\Q)$ is
exact. Thus, the reduced group in which $c^+_{\red}(\xi;[\omega])$
lives consists only of torsion classes, and the result now follows
from Theorem~\ref{thm:HFaWeakFilling}.  \endproof

Sometimes, it is easier to use $\Zmod{p}$ coefficients (especially
when $p=2$). To this end, we say that $Y$ a rational homology
three-sphere is a $\Zmod{p}$--$L$--space for some prime $p$ if
$\HFa(Y;\Zmod{p})$ has rank $|H_1(Y;\Z)|$ over $\Zmod{p}$ (of course,
an $L$ space is automatically a $\Zmod{p}$--$L$--space for all
$p$). Since $c^+(\xi;[\omega])$ is primitive, the above argument shows
that a $\Zmod{p}$--$L$--space (for any prime $p$) cannot support a taut
foliation.

The need to use twisted coefficients in the statement of
Theorem~\ref{thm:HFaWeakFilling} is illustrated by the three-manifold
$Y$ obtained as zero-surgery on the trefoil. The reduced Heegaard
Floer homology with untwisted coefficients is trivial
(cf.\ Equation~26
of~\cite{AbsGraded}), but this three-manifold admits a taut foliation. 
(In particular the  reduced Heegaard Floer homology of this manifold
with twisted
coefficients is non-trivial, cf.\ Lemma~8.6
of~\cite{AbsGraded}.)

\section{The Thurston norm}
\label{sec:ThurstonNorm}

We turn our attention to the proof of Theorem~\ref{thm:ThurstonNorm}.

\medskip
\noindent{\bf{Proof of Theorem~\ref{thm:ThurstonNorm}}\qua}
It is shown in Section~1.6
of~\cite{HolDiskTwo} that if $\uHFa(Y,\spinc)\neq 0$, then 
\begin{equation}
\label{eq:AdjunctionInequality}
|\langle
c_1(\spinc),\xi\rangle | \leq \Theta(\xi).
\end{equation} 
(The result is stated there for $\HFp$ with untwisted coefficients,
but the argument there applies to the case of $\uHFa$.)  It remains to
prove that if $\Sigma\subset Y$ is an embedded surface which
minimizes complexity in its homology class $\xi$, then there is a $\SpinC$
structure $\spinc$ with $\uHFa(Y,\spinc)\neq 0$ and 
\begin{equation}
\label{eq:AdjunctionSharp}
\langle
c_1(\spinc),[\Sigma]\rangle = -\chi_+(\Sigma).
\end{equation}

The K\"unneth principle for connected sums
(cf.\ Theorem~1.5 of~\cite{HolDiskTwo})
states that $$\uHFa(Y_1\# Y_2,\spinc_1\# \spinc_2)\otimes_\Z \Q
\cong 
\uHFa(Y_1,\spinc_1)\otimes_\Z \uHFa(Y_2,\spinc_2)\otimes_\Z\Q.$$
In particular, if $\uHFa(Y_1,\spinc_1)\otimes_{\Z}\Q$ and
$\uHFa(Y_2,\spinc_2)\otimes_{\Z}\Q$ are non-trivial, then so is
$\uHFa(Y_1\# Y_2,\spinc_1\#\spinc_2)\otimes_\Z \Q$. Since every closed
three-manifold admits a connected sum decomposition where the summands
are all either irreducible or copies of $S^2\times S^1$~\cite{Milnor},
it suffices to verify that $\uHFa(Y,\spinc)\otimes_\Z \Q$ is
non-trivial for the elementary summands of $Y$. (It is straightforward
to see that $\Theta_{Y_1\# Y_2}(\xi_1+\xi_2) =
\Theta_{Y_1}(\xi_1)+\Theta_{Y_2}(\xi_2)$ in $Y_1\#Y_2$, where here
$\xi_i\in H_2(Y_i)$, under the natural identification $H_2(Y_1\#
Y_2)\cong H_2(Y_1)\oplus H_2(Y_2)$.)

We first observe that if $Y$ has trivial Thurston semi-norm (for
example, when $b_1(Y)=0$ or $Y=S^2\times S^1$), then there is an
element $\spinc\in\SpinC(Y)$ for
which $\uHFa(Y,\spinc)\neq 0$. Indeed, it is shown in
Theorem~10.1 of~\cite{HolDiskTwo} that
$\uHFinf(Y,\spinc)\cong \Z[U,U^{-1}]$ for any $\spinc$ with
$c_1(\spinc)=0$. Also, for such $\SpinC$ structures, the map from
$\uHFinf(Y,\spinc)$ to $\uHFp(Y,\spinc)$ is non-trivial.  The
non-triviality of $\uHFa(Y,\spinc)$ follows at once (using the
analogue of Exact Sequence~\eqref{eq:ExactSequences} for the case of
twisted coefficients).

In the case where $Y$ is an irreducible three-manifold
with non-trivial Thurston norm, 
and $\Sigma$ is a surface which minimizes
complexity in its homology class, Gabai~\cite{Gabai} constructs a
smooth taut foliation ${\mathcal F}$ for which $$\langle c_1({\mathcal
F}),[\Sigma]\rangle = -\chi_+(\Sigma).$$ 
According to a theorem of
Eliashberg and Thurston, then $[-1,1]\times Y$ can be equipped with a
convex symplectic form, which extends ${\mathcal F}$, thought of as a
foliation over $\{0\}\times Y$. In particular, their result gives a
weakly symplectically semi-fillable contact structure $\xi$ with
$\langle c_1(\xi),[\Sigma]\rangle = -\chi_+(\Sigma)$. It follows now from
Theorem~\ref{thm:HFaWeakFilling} that $c(\xi,[\omega])\in
\uHFa(Y,[\omega],\spinc(\xi))\otimes_\Z \Q\neq 0$.
\endproof

One approach to Theorem~\ref{thm:Knots} would directly relate knot
Floer homology with the twisted Floer homology of the zero-surgery. We
opt, however, to give an alternate proof which uses the relation
between the knot Floer homology and the Floer homology of the
zero-surgery in the untwisted case, and adapts the proof rather than the
statement of Theorem~\ref{thm:ThurstonNorm}.  The 
relevant relationship between these groups can be
found in Corollary~4.5
of~\cite{HolDiskKnots}, according to which if
$d>1$ is the smallest integer for which $\HFKa(K,d)\neq 0$, then
\begin{equation}
\label{eq:RelateZeroSurgery}
\HFKa(K,d)\cong \HFp(S^3_0(K),d-1),
\end{equation} where here we have identified
$\SpinC(S^3_0(K))\cong \Z$ by the map $\spinc \mapsto \langle
c_1(\spinc),[\Sigma]\rangle/2$, where  $[\Sigma]\in
H_2(S^3_0(K);\Z)\cong \Z$ is some generator.
(Note that the choice of generator is not particularly important, as
$\HFp(S^3_0(K),i)\cong \HFp(S^3_0(K),-i)$, according to 
the conjugation invariance of Heegaard Floer homology,
Theorem~2.4 of~\cite{HolDiskTwo}.)

This result will be used in conjunction with the ``adjunction inequality''
for knot Floer homology,  Theorem~5 of~\cite{HolDiskKnots},
which shows that $\HFKa(K,i)=0$ for all $|i|>g(K)$; and indeed, the
proof of that result proceeds by constructing a compatible doubly-pointed Heegaard
diagram (from a genus-minimizing Seifert surface for $K$) which has no
simultaneous trajectories $\x$ with $s(\x)>g(K)$. 

\medskip
\noindent{\bf{Proof of Theorem~\ref{thm:Knots}}}\qua
Let $K\subset S^3$ be a knot with genus $g$. Assume for the moment
that $g>1$. Let $Y$ be the three-manifold obtained as zero-framed
surgery on $S^3$ along $K$, and let $[\Sigma]\in H_2(Y;\Z)$ denote a
generator.  In this case, Gabai~\cite{GabaiKnots} constructs a taut
foliation ${\mathcal F}$ over $Y$ with $\langle c_1({\mathcal F}),
[\Sigma]\rangle = 2-2g$. Eliashberg's theorem~\cite{Eliashberg} now
provides a symplectic four-manifold $X=X_1\cup_Y X_2$, where here
$b_2^+(X_i)>0$. According to the product formula
Equation~\eqref{eq:ProductFormula}, the sum $$\sum_{\eta\in H^1(Y)}
\Phi_{X,\spinccan(\omega)+\delta\eta}$$ is calculated by a
homomorphism which factors through the Floer homology
$\HFp(Y,\spinccan(\omega)|_{Y})$. On the other hand,
$c_1(\spinccan(\omega))$ gives a cohomology class whose evaluation on
a generator for $H_2(Y;\Z)$ is non-trivial when $g>1$ (for a suitable
generator, this evaluation is given by $2-2g$). Since the image of 
a generator of $H^1(Y;\Z)$ is represented by a surface in $X$ with 
square zero and non-zero evaluation of $c_1(\spinc(\omega))$, it follows
that the various terms in the sum are homogeneous of different degrees.  
But by 
Theorem~1.1 of~\cite{HolDiskSymp},
it follows that the term corresponding to $\spinccan(\omega)$
(and hence the sum) is non-trivial. It
follows now that $\HFp(Y,\spinccan(\omega)|_{Y})=\HFp(S^3_0(K),g-1)$
(for suitably chosen generator) is non-trivial and hence, in view of
Equation~\eqref{eq:RelateZeroSurgery}, Theorem~\ref{thm:Knots}
follows for knots with genus at least two.

Suppose that $g=1$. In this case, we have a K\"unneth principle for
the knot Floer homology (cf.\ Equation~5
of~\cite{HolDiskKnots}), according to which (since $\HFKa(K,s)=0$ for all
$s>1$), $$\HFKa(K\# K,2)\otimes_\Z\Q\cong \HFKa(K,1)\otimes_\Q \HFKa(K,1).$$ But
$K\#K$ is a knot with genus $2$, and hence $\HFKa(K\# K,2)$ is
non-trivial; and hence, so is $\HFKa(K,1)$.
\endproof

\medskip
\noindent{\bf{Proof of Corollary~\ref{cor:KMOSz}}}\qua
According to the integral surgeries long exact sequence for Heegaard
Floer homology (in its graded form), if $S^3_p(K)\cong L(p,1)$, the
Alexander polynomial of $K$ is trivial (indeed
$\HFp(S^3_0(K))\cong\HFp(S^2\times S^1)$),
cf.\ Theorem~1.8
of~\cite{AbsGraded}. In~\cite{NoteLens}, it is shown that if
$S^3_p(K)$ is a lens space for some integer $p$, then the knot Floer
homology $\HFKa_*(K,*)$ is determined by the Alexander polynomial
$\Delta_K(T)$ (cf.\ Theorem~1.2
of~\cite{NoteLens}) which in the present case is trivial. Thus, in
view of Theorem~\ref{thm:Knots}, the knot $K$ is trivial.
\endproof

\medskip
{\bf{Proof of Corollary~\ref{cor:MorseInterp}}}\qua
In the proof of Theorem~5 of~\cite{HolDiskKnots}, we
demonstrate that if a knot has genus $g$, then there is a compatible
Heegaard diagram with no simultaneous trajectories $\x$ for which 
$s(\x)>g$. In the opposite direction, note that $\HFKa(K,d)$ is generated
by simultaneous trajectories with $s(\x)=d$. According to 
Theorem~\ref{thm:Knots}, $\HFKa(K,g)\neq 0$, and hence any
compatible Heegaard diagram must contain some simultaneous trajectories $\x$
with $s(\x)=g$.
\endproof

\end{document}

%% file: gtoutput.tex
\def\ifplaintex{\expandafter\ifx\csname documentclass\endcsname\relax}

\ifplaintex 
\else
\fi

\expandafter\ifx\csname beginpicture\endcsname\relax
\expandafter\ifx\csname documentclass\endcsname\relax
\input pictex \else
\input prepictex \input pictex \input postpictex \fi\fi

\def\gt{{\mathsurround=0pt\it $\cal G\mskip-2mu$eometry \&\ 
$\cal T\!\!$opology}}        

\def\gtp{{\mathsurround=0pt\it $\cal G\mskip-2mu$eometry \&\ 
$\cal T\!\!$opology $\cal P\!$ublications}}  

\def\lognumber#1{\def\thelognumber{#1}}
\def\volumenumber#1{\def\thevolumenumber{#1}}
\def\papernumber#1{\def\thepapernumber{#1}}
\def\volumeyear#1{\def\thevolumeyear{#1}}

\def\pagenumbers#1#2{\def\startpage{#1}\def\finishpage{#2}}
\def\published#1{\def\publishdate{#1}}
\def\proposed#1{\def\theproposer{#1}}
\def\seconded#1{\def\theseconders{#1}}
\def\received#1{\def\receiveddate{#1}}
\def\revised#1{\def\reviseddate{#1}}
\def\accepted#1{\def\accepteddate{#1}}

\def\coverauthors#1{\def\thecoverauthors{#1}}
\def\asciiauthors#1{\def\theasciiauthors{#1}}
\def\asciiaddress#1{\def\theasciiaddress{#1}}
\def\asciiemail#1{\def\theasciiemail{#1}}

\long\def\asciiabstract#1{\long\def\theasciiabstract{#1}}

\let\\\par\let\thelognumber\relax
\let\thevolumenumber\relax\let\thepapernumber\relax
\let\thevolumeyear\relax\let\thesamplenumber\relax\let\startpage\relax
\let\finishpage\relax\let\publishdate\relax\let\receiveddate\relax
\let\reviseddate\relax\let\accepteddate\relax\let\theasciititle\relax
\let\theasciiauthors\relax\let\theasciiaddress\relax
\let\theasciiabstract\relax
\let\theasciiemail\relax\let\theshortauthors\relax\let\theshorttitle\relax
\let\thecoverauthors\relax

\long\def\maketitlep{   

\count0=\startpage

\gt\hfill      
\beginpicture
\setcoordinatesystem units <0.33truein, 0.33truein> point at 2.2 0.9
\setplotsymbol ({$\cal G$})
\plotsymbolspacing=9truept
\circulararc 315 degrees from 0 1 center at 0 0
\setplotsymbol ({$\cal T$})
\circulararc 315 degrees from 1 -1 center at 1 0
\endpicture
\break
{\small\ifx\thesamplenumber\relax 
Volume \else Sample
\fi\thevolumenumber\ (\thevolumeyear)
\startpage--\finishpage\nl
Published: \publishdate}
\vglue 0.5truein plus 0.4fil minus 0.1truein

{\parskip=0pt\leftskip 0pt plus 1fil\def\\{\par\smallskip}{\ifplaintex\large
\else\Large\fi\bf\thetitle}\par\medskip}   

\vglue 0pt plus 0.1fil 

{\parskip=0pt\leftskip 0pt plus 1fil\def\\{\par}{\sc\theauthors}
\par\medskip}

\vglue 0pt plus 0.1fil 

{\small\parskip=0pt\let\newline\\
{\leftskip 0pt plus 1fil\def\\{\par}{\sl\theaddress}\par}
\expandafter\ifx\theemail\relax    
\relax\else\vglue 5pt plus 0.02fil minus 2pt\def\\{\stdspace{\rm 
and}\stdspace} 
\cl{Email:\stdspace\tt\theemail}\fi
\ifx\theurl\relax                  
\relax\else\vglue 5pt plus 0.02fil minus 2pt\def\\{\stdspace{\rm 
and}\stdspace}
\cl{URL:\stdspace\tt\theurl}\fi\par}

\vglue 7pt plus 0.3fil minus 3pt

{\bf Abstract}
\vglue 5pt plus 0.1fil minus 2pt

\theabstract

\vglue 7pt plus 0.3fil minus 3pt

{\bf AMS Classification numbers}\quad Primary:\quad \theprimaryclass

Secondary:\quad \thesecondaryclass

\vglue 5pt plus 0.3fil minus 2pt

{\bf Keywords:}\quad \thekeywords

\vglue 10pt plus 0.5fil minus 5pt

{\small  Proposed: \theproposer\hfill Received: \receiveddate\nl
Seconded: \theseconders\hfill 
\ifx\reviseddate\relax                         
Accepted: \accepteddate                        
\else
Revised: \reviseddate                          
\fi}
\eject
}       

\let\maketitlepage\maketitlep
\let\maketitle\maketitlepage

\font\phead=cmsl9 scaled 950
\font=cmsl9 scaled 1050
\font\pnum=cmbx10 scaled 913
\font\lnum=cmbx10 
\font\pfoot=cmsl9 scaled 950
\font=cmsl9 scaled 1050
\ifplaintex
\headline{\vbox to 0pt{\vskip -4.5mm\line{\small\phead\ifnum
\count0=\startpage ISSN 1364-0380 (on line)
1465-3060 (printed) \hfill {\pnum\folio}\else\ifodd\count0\def\\{ }%
\ifx\theshorttitle\relax\thetitle\else\theshorttitle\fi\hfill{\pnum\folio}
\else\def\\{ and }{\pnum\folio}\hfill\ifx\theshortauthors\relax\theauthors
\else\theshortauthors\fi\fi\fi}\vss}}
\footline{\vbox to 0pt{\vglue 0mm\line{\small\pfoot\ifnum\count0=\startpage
\copyright\ \gtp\hfill\else
\gt, Volume \thevolumenumber\ (\thevolumeyear)\hfill\fi}\vss
}}
\else
\makeatletter
\def\@oddhead{{\small\lhead\ifnum\count0=\startpage ISSN 1364-0380 (on line)
1465-3060 (printed) \hfill {\lnum\number\count0}\else\ifodd\count0
\def\\{ }\ifx\theshorttitle\relax \thetitle \else\theshorttitle\fi\hfill
{\lnum\number\count0}\else\def\\{ and }{\lnum\number\count0}
\hfill\ifx\theshortauthors\relax 
\theauthors\else\theshortauthors\fi\fi\fi}}\def\@evenhead{\@oddhead}
\def\@oddfoot{\small\lfoot\ifnum\count0=\startpage\copyright\ \gtp\hfill\else
\gt, Volume \thevolumenumber\ (\thevolumeyear)\hfill\fi}
\def\@evenfoot{\@oddfoot}
\makeatother
\fi

\newwrite\gtoutfile
\long\gdef\makeheadfile{  
{\def\\{, }\def\s{ }
\immediate\openout\gtoutfile head.xxx
\immediate\write\gtoutfile{Proxy-for: \ifx\theasciiauthors\relax
\theauthors\else\theasciiauthors\fi\s<\ifx\theasciiemail\relax\theemail\else\theasciiemail\fi>}
\immediate\write\gtoutfile{\noexpand\\}
\immediate\write\gtoutfile{Authors: \ifx\theasciiauthors\relax
\theauthors\else\theasciiauthors\fi}
{\def\\{ }\immediate\write\gtoutfile{Title: \ifx\theasciititle\relax
\thetitle\else\theasciititle\fi}}
\immediate\write\gtoutfile{Subj-class: GT or SG or MG etc}
\immediate\write\gtoutfile{MSC-class: \theprimaryclass\ifx\thesecondaryclass\relax\else, \thesecondaryclass\fi}
\immediate\write\gtoutfile{Journal-ref: Geom. Topol. \thevolumenumber
(\thevolumeyear) \startpage-\finishpage}
\immediate\write\gtoutfile{Comments: Published by Geometry and Topology at}
\immediate\write\gtoutfile{\s\s http://www.maths.warwick.ac.uk/gt/GTVol\thevolumenumber/paper\thepapernumber.abs.html}
\immediate\write\gtoutfile{\noexpand\\}
\immediate\write\gtoutfile{}
\ifx\theasciiabstract\relax
\immediate\write\gtoutfile{\theabstract}\else
\immediate\write\gtoutfile{\theasciiabstract}\fi
\immediate\write\gtoutfile{}
\immediate\write\gtoutfile{\noexpand\\}
\immediate\write\gtoutfile{}
\immediate\closeout\gtoutfile}}  

\def\maketitlepage{\maketitlep\makeheadfile}
\let\maketitle\maketitlepage